\documentclass [12pt]{article}

\begin{document}
\title{LEAST SQAURE METHOD FOR SUM OF THE FUNCTIONS SATYSFYING THE DIFFERENTIAL EQUATIONS WITH POLYNOMIAL COEFFICIENTS.}

\author{Berngardt O.I., Voronov A.L.}

\maketitle

\begin{abstract}
We propose a linear algorithm for determining two function 
parameters by their linear combination. These functions must 
satisfy the first order differential equations with polynomial 
coefficients and our parameters are the coefficients of these 
polynomials. The algorithm consists of sequential solution by 
least squares method of two linear problems - first, differential 
equation polynomial coefficients  determining for linear combination 
of two given functions and second - determining functions parameters 
by these polynomial coefficients. Numerical modeling carried by 
this scheme gives an good accordance under weak normal noise (with 
dispersion \( <5\% \)).
\end{abstract}

\textbf{Introduction}

In practice we often run against a problem of determining two 
functions parameters from certain classes by their linear 
combination. 
In physics problems of that kind usually called the problems of signals, modes or rays separation. The examples of these certain classes are exponents with polynomial arguments, harmonics functions with polynomial phase and rational functions \cite{6}. Such problems can arise from experiment interpretation when function class in linear combination is derived from some theoretical consideration , but parameters of these functions are unknown.

One of the most common methods for functions parameters determining is the least squares method, which in general case boils down to absolute minimization of the function from given class which depends on many unknown parameters. In linear case these problem has a simple solution - essentially, the linear system equation solution for given matrix. In our case, when parameters dependence is nonlinear, this problem is sufficiently hard and does not have an analytical solution and so (when direct parameters search is not acceptable by resources conditions) the problem usually attacks with some step-by-step function quasilinearization (Newton's methods)\cite{1}. However, these methods has some serious faults - strong dependence on initial approximation and need for big computational resources. These conditions cannot be met if each function computing needs too much time or if we do not have good initial approximation and we are not sure in unique global minimum existence.

We propose an algorithm for transformation of the above-mentioned problem to problem of finding minimum for functionals with linear dependence on parameters and so we may apply the well-known methods based on linear systems equations solution. We restrict our attention to the case of functions  satisfying the first order differential equation with polynomial coefficients. This problem becomes linear if we shifted from function parameters determining to the problem  of  differential equations parameters determining (due to the fact that differential equations are linear with respect to these polynomial coefficients) 

The method similar to ours was developed by Kulikov \cite{2} for differential equations with constant coefficients,
but in our case this method cannot be applied because our function does not satisfy differential equation with constant coefficients.

\textbf{1. The problem.}

Let \(f_1 (t),f_2 (t)\) - unknown vector-valued functions satisfying the differential equation with polynomial coefficients :

\begin{equation}
\label{eq:1}
[D + \hat P_i (t)]f_i (t) = 0
\end{equation}

where

\begin{equation}
\label{eq:2}
\hat P_i (t) = \hat M_i \sum\limits_{j = 0}^{N_i } {p_{ij} t^j } 
\end{equation}

- polynomial coefficients of degree \( N_1  \ge N_2 \)  with unknown coefficients 
\( p_{ij} \), matrixes \( \hat M_1 ,\hat M_2 \) commute (\( \hat M_1 \hat M_2  = \hat M_2 \hat M_1 \) ) and are given, but operators themselves does not equal , \( \hat P_1 (t) \ne \hat P_2 (t) \). Suppose that we know the linear combination of these functions on interval \( t \in [0,1] \):

\begin{equation}
\label{eq:3}
F(t) = a_1 f_1 (t) + a_2 f_2 (t),
\end{equation}

and that linear combination is bounded on this interval (\( \sup (F(t)) < \infty  \)), but coefficients \( a_1 ,a_2  \ne 0 \) are unknown.

Our goal is to transform the problem of determining nonlinear parameters \( p_{ij} \) in model function 

\begin{equation}
\label{eq:4}
F_{\bmod } (t;p_{ij} ) = a_1 f_1 (t;p_{1j} ) + a_2 f_2 (t;p_{2j} )
\end{equation}

to the sequential series of linear problems. 

\textbf{2. The algorithm description.}

\textbf{Main idea.} Consider the problem of determining \( N_1  + N_2  + 2 \)
unknown parameters \( p_{ij} \) for two vector-valued functions by their linear combination (3). From (1)-(2) we can see that unknown parameters \( p_{ij} \) are nonlinear one for functions \( f_i (t;p_{ij} ) \). So we have a problem of determining of \( N_1  + N_2  + 2 \) unknown parameters \( p_{ij},a_{1},a_{2} \) by given function \( F(t) \), and only two of these parameters are linear for this function. 

One of most common methods for determining  the unknown parameters by given function is the least squares method , which in our case essentially the minimization of deviation functional by all set of parameters 

\begin{equation}
\label{eq:5}
\Omega (p_{ij} ,a_1 ,a_2) = \int\limits_0^1 {(F(t) - F_{\bmod } (t;p_{ij} ,a_1 ,a_2 ))^2 dt} 
\end{equation}

where \( F _{\bmod } (t;p_{ij},a_{1},a_{2} ) \) - model function (4) with given parameters dependence and \( F(t) \) - experimental function.

Least squares method gives simple analytical solution which stable small errors in \( F(t) \)  when all the parameters for functional  (\ref{eq:5})  are linear for model function (\ref{eq:4}). In our case that is not so and thus direct least squares approach is not acceptable because of complex functional behavior and need for  big computational resources. So we come to idea to develop a linear algorithm due to specific properties of our problem. 

Although our parameters \( p_{ij} ,a_{1} ,a_{2} \) are nonlinear for model function (\ref{eq:4}), differential equations (\ref{eq:1}) are linear with respect to \( p_{ij} \)  and does not depend on scalar factor for \( f_{i }(t) \). Thus it makes sense to analyze differential equations themselves and this approach may be more constructive than direct approach. 

Formally, the problem of two functions parameters determining boils down to  determining of the unknown coefficients  \( p_{ij} \) in differential equations (\ref{eq:1}) by given linear combination of unknown functions (\ref{eq:3}). If the function (\ref{eq:3}) satisfies the differential equation \( \hat L(t,p_{ij} )F_{\bmod } (t;p_{ij} ) = 0 \), where all the unknown parameters are linear, then we may obtain these parameters by minimization of the functional:

\begin{equation}
\label{eq:6}
\Omega (p_{ij} ) = \int\limits_0^1 {(\hat L(t,p_{ij} )F(t))^2 dt} 
\end{equation}

In order to diminish the influence of numerical differentials we may integrate our differential equation thus transform it to integral equation. By solving the linear problem for functional (\ref{eq:6}) we obtain the unknown parameters   for functions \( f_{i }(t) \). If these parameters are obtained, then our initial problem (\ref{eq:3}) is linear for remaining parameters \( a_{1} ,a_{2} \) and so it has an analytical solution. 

Differential equation for sum of unknown functions. Our next goal is to construct a differential equation for function \( F(t) \). Our initial equation are homogeneous, and so differential equation for sum (\ref{eq:3}) does not depend on parameters \( a_{1} ,a_{2} \).

Suppose that we have an equation: 

\begin{equation}
\label{eq:7}
[\hat AD^2  + \hat BD + \hat C]F = 0
\end{equation}

Then we have:

\begin{equation}
\label{eq:8}
[\hat AD^2  + \hat BD + \hat C]f_{i} = 0
\end{equation}

From system (\ref{eq:8}) we obtain a linear system for unknown functions \( A(t),B(t),C(t) \). If the operators \( \hat P_i  \) commute, then this system has a simple solution:

\begin{equation}
\label{eq:9}
\left\{ \begin{array}{l}
 A = P_2  - P_1  \\ 
 B = (P_2^2  - P_1^2 ) - \frac{d}{{dt}}(P_2  - P_1 ) \\ 
 C = P_2 P_1 (P_2  - P_1 ) - (P_1 \frac{{dP_2 }}{{dt}} - P_2 \frac{{dP_1 }}{{dt}}) \\ 
 \end{array} \right.
\end{equation}

As far as our functions \( P_1 (t),P_2 (t) \) are polynomials of given degrees, then the functions (\ref{eq:9})  also are polynomials; for example: 

\begin{equation}
B_k (t) = \sum\limits_{k = 0}^{N_1  + N_2 } {\beta _k t^k \hat Y_k } 
\end{equation}

where \( \beta _k \) - new coefficients, which has explicit dependence on parameters \( p_{ij} \), and \( \hat Y_k \) - matrixes which are explicitly obtained from matrixes 
\( \hat M_1 ,\hat M_2 \). 

By using differential equation  (\ref{eq:7}),(\ref{eq:9}) for linear combination of initial functions (\ref{eq:3}) we may construct two-step algorithm for determining parameters , with each step linear. In order to do this, transform our problem to two linear problems. 

\textbf{First step.} Our first goal is to determine the coefficients of the polynomials \( A(t),B(t),C(t)\) under assumption that all these coefficients are independent. This problem is linear since all the coefficients are linear in differential equation (\ref{eq:7}):

\begin{equation}
\label{eq:7.1}
[\sum\limits_{k = 0}^{N_1 } {\alpha _k t^k \hat X_k } D^2  + \sum\limits_{k = 0}^{N_1  + N_2 } {\beta _k t^k \hat Y_k } D + \sum\limits_{k = 0}^{2N_1  + N_2 } {\gamma _k t^k \hat Z_k } ]F = 0
\end{equation}

So the problem of obtaining the coefficients \( \alpha _k ,\beta _k ,\gamma _k  \) is the linear least squares problem for \( 4N_{1}+2N_2+3\) linear parameters, and we easily solve this problem \cite{4} . 

In order to diminish the numerical errors let us transform our differential equation (\ref{eq:7}) to integral one by integrating it two times: 

\begin{equation}
\label{eq:10}
\begin{array}{l}
 \hat A(t)F(t) + \int\limits_0^t {\left( {\hat B(x) - 2\frac{{d\hat A(x)}}{{dx}}} \right)F(x)dx}  +  \\ 
  + \int\limits_0^t {dx\int\limits_0^x {dy} \left( {\hat C(y) - \frac{{d\hat B(y)}}{{dy}} + \frac{{d^2 \hat A(y)}}{{d^2 y}}} \right)F(y)dy + \hat L_1 t = \hat L_0 }  \\ 
 \end{array}
\end{equation}

This integral equation also linear with respect to parameters \( \alpha _k ,\beta _k ,\gamma _k  \), and all its coefficients are known functions - certain combinations of the function \( F(t)\) and its integrals. Since all the coefficients are known and additional integration constants (which are also must be determined) are linear, then linear least squares method gives us an analytical solution. To exclude the trivial solution we must set one of the parameters (for example \( L_{0}\)) equal to 1 and after that we may find parameters \( \alpha _k ,\beta _k ,\gamma _k , L_{1} \) by minimization of the functional:

\begin{equation}
\label{eq:11}
\begin{array}{l}
 \Omega (L_1 ,\alpha _k ,\beta _k ,\gamma _k ) = \int\limits_0^1 {dt} \left( {\hat A(t)F(t) + \int\limits_0^t {\left( {\hat B(x) - 2\frac{{d\hat A(x)}}{{dx}}} \right)F(x)dx}  + } \right. \\ 
  + \left. {\int\limits_0^t {dx\int\limits_0^x {dy} \left( {\hat C(y) - \frac{{d\hat B(y)}}{{dy}} + \frac{{d^2 \hat A(y)}}{{d^2 y}}} \right)F(y)dy}  + \hat L_1 t - 1} \right)^2  = \min  \\ 
 \end{array}
\end{equation}

Since this problem is linear with respect to all parameters \( \alpha _k ,\beta _k ,\gamma _k , L _{1} \), we may obtain the solution by linear least squares method. 

\textbf{Second step.} Our second goal is to determine the unknown coefficients of the initial polynomials  \( \hat P_1 (t),\hat P_2 (t) \) by given functions \( A(t),B(t),C(t) \). This problem also ca transform to the linear problem. Really, after solving the problem (\ref{eq:11}) for functions \( A(t),B(t),C(t) \), we have three functions: 

\begin{equation}
\label{eq:12}
\begin{array}{l}
 K_1 (t) = A(t) = P_2  - P_1  \\ 
 K_2 (t) = B(t) + \frac{{dA(t)}}{{dt}} = P_2^2  - P_1^2  \\ 
 K_3 (t) = C(t) = P_2 P_1 (P_2  - P_1 ) - (P_1 \frac{{dP_2 }}{{dt}} - P_2 \frac{{dP_1 }}{{dt}}) \\ 
 \end{array}
\end{equation}

Under our assumptions we obtain functions \( K_1 (t),K_2 (t),K_3 (t) \) only up to constant factor (due to excluding of the trivial solution). So we derive from system (\ref{eq:12}) two equations which is invariant under multiplication of the functions \( K_1 (t),K_2 (t),K_3 (t) \) by the constant factor: 

\begin{equation}
\label{eq:13}
\begin{array}{l}
 P_1 (t) + P_2 (t) = K_2 /K_1  \\ 
 P_1 (t)P_2 (t) = [K_3  - (K_1 \frac{d}{{dt}}(K_2 /K_1 ) - (K_2 /K_1 )\frac{d}{{dt}}K_1 )]/K_1  \\ 
 \end{array}
\end{equation}

From this system we can easily obtain an analytical expressions for functions
\( \hat P_1 (t),\hat P_2 (t) \):  

\begin{equation}
\label{eq:14}
P_{1,2} (t) = S_{1,2} (t) 
\end{equation}

where \( S_{1,2} (t) \) has explicit expression as functions of \( K_1 (t),K_2 (t),K_3 (t) \) due to Viet theorem (\ref{eq:13}):

\begin{equation}
\label{eq:15}
\begin{array}{l}
 S_1 (t) + S_2 (t) = K_2 /K_1  \\ 
 S_1 (t)S_2 (t) = [K_3  - (K_1 \frac{d}{{dt}}(K_2 /K_1 ) - (K_2 /K_1 )\frac{d}{{dt}}K_1 )]/K_1  \\ 
 \end{array}
\end{equation}

The coefficients of these polynomials can be found from functional minimum condition:

\begin{equation}
\label{eq:16}
 \Omega _i (p_{ij} ) = \int\limits_0^1 {(P_i (t;p_{ij} ) - S_i (t))^2 dt}, 
\end{equation}

where analytical expression for \( P_i (t;p_{ij} ) \) are given by (\ref{eq:2}), and \( S_i (t) \) are determined by the first step (\ref{eq:15}).

Since our model polynomials \( P_i (t;p_{ij} ) \) are linear with respect to \( p_{ij} \) (\ref{eq:2}), and \( S_i (t) \)  are determined by the first step (\ref{eq:15}), the problem of determining the polynomial coefficient is linear and so can be solved by linear least squares method \cite{4}. Thus our algorithm transforms initial nonlinear problem to sequential two-step linear problem for differential equation (\ref{eq:7}) and each step is linear least squares problem for certain set of parameters. 

\textbf{3. Modeling results.}

Our algorithm have been tested for two problems - separation of two overlapping Gaussian functions \( f_i (t) = R_i \exp ( - (\alpha _i t^2  + \beta _i t)) \), satisfying the differential equation 

\begin{equation}
[D - (2\alpha _i t + \beta _i )]f_i  = 0,
\end{equation}

and separation of two vector-valued signals with linear frequency modulation 

\begin{equation}
f_i (t) = R_i \left( \begin{array}{l}
 \cos (\alpha _i t + \beta _i t^2 ) \\ 
 \sin (\alpha _i t + \beta _i t^2 ) \\ 
 \end{array} \right) ,
\end{equation}

satisfying the differential equation
\begin{equation}
\left[ {D + (\alpha _i  + 2\beta _i t)\left( {\begin{array}{*{20}c}
   0 & { - 1}  \\
   1 & 0  \\
\end{array}} \right)} \right]f_i  = 0
\end{equation}

In the first case algorithm is valid under \( 10-15\% \) normal noise (For function maximal amplitude), in the second case - under \(3-5\%\) normal noise.

\textbf{4. Discussion.}

In practice we often have the case when operators in equation (\ref{eq:1}) are not polynomials but a rational functions:

\begin{equation}
\hat P_i (t) = \hat M_i \frac{{\sum\limits_{j = 0}^{N_{pi} } {p_{ij} t^j } }}{{\sum\limits_{j = 0}^{N_{qi} } {q_{ij} t^j } }} ,
\end{equation}

under assumption that \( P_i (t) \) has no singularities.

Such class includes rational functions, exponents and harmonic functions with rational argument. Our method also can be applied to that extended class. The reason for this is that functions \( K_1 (t),K_2 (t),K_3 (t) \) will be also rational so after multiplying our differential equation (\ref{eq:7}) by some polynomial we can obtain the differential equation with polynomial coefficients . Since determining \( P_i (t) \) by \( K_1 (t),K_2 (t),K_3 (t) \) does not depend on constant factor (even if this factor is nonconstant polynomial), then the expressions (\ref{eq:14}-\ref{eq:15}) are also valid in this For determining all the coefficients we must multiply all the expressions by \( P_i (t) \) denominator, and set one of the coefficients \( q_{ij} \) (for example, under highest degree) to 1 to exclude the trivial solution:

\begin{equation}
\hat M_i \sum\limits_{j = 0}^{N_{pi} } {p_{ij} t^j }  - \sum\limits_{j = 0}^{N_{qi}  - 1} {q_{ij} t^j } \hat S_i (t) = t^{N_{qi} } \hat S_i (t)
\end{equation}

This problem also linear with respect to parameters \( p_{ij} , q_{ij}\) and thus can be solved by linear least squares method.

\textbf{Conclusion.}

We propose an algorithm for determining parameters of two functions satisfying differential equations of the first order with polynomial coefficients by their linear combination. This algorithm transform initial nonlinear problem to sequential solving of the two linear problems - determining of polynomial coefficients for linear differential equation of the second order with polynomial coefficients, and determining initial polynomial coefficients from results or the first step. Numerical testing shows us that the algorithm is valid under weak normal noise  (with dispersion \(<5\% \)).

\textbf{Acknowledgements.}

Thanks are due to V.E.Nosov for interest to our work and to I.I.Orlov for fruitful discussions. This work was partially supported by RBFR grant 'Fundamental science schools support' No 00-15-98509.


\begin{thebibliography}{10}
\bibitem{1}
J.E.Dennis, R.B.Schnabel Numerical methods for unconstrained optimization and nonlinear equations, Prentice-Hall Inc., 1983, 430p.
\bibitem{2}
N.K.Kulikov, G.N. Bagautdinov, Normal differential equations. The solution of differential equations based on the functions with flexible structure(in Russian)- Alma-Ata, 1973, 132p. 
\bibitem{4}
J.H.Wilkinson, C.Reinsch, Handbook for Automatic Computation Linear Algebra, Heidelberg, 1969, 360pp.
\bibitem{5}
G.A.Korn, T.M.Korn, Mathematical handbook for scientists and enginiers, McGraw-Hill Book Company, 1968, 830p.
\bibitem{6}
6. E. Von Kamke, Handbook of normal differential equations, (in Russian) Moscow, Nauka, 1976, 576pp. 
\end{thebibliography}
\end{document}